\documentclass{autart}

\usepackage{amssymb,amsmath}
\usepackage{natbib}

\def\fudge{\mathchoice{}{}{\mkern.5mu}{\mkern.8mu}}
\def\bbc#1{{\rm l\mkern-9.4mu \mbox{#1}}}
\def\bbb#1{{\rm I\mkern-3.5mu #1}}
\def\bba#1#2{{\rm #1\mkern-#2mu\fudge #1}}
\def\bb#1{{\count4=`#1 \advance\count4by-64 \ifcase\count4\or\bba A{11.5}\or
   \bbb B\or\bbc C\or\bbb D\or\bbb E\or\bbb F \or\bbc G{5}\or\bbb H\or
   \bbb I\or\bbc J{3}\or\bbb K\or\bbb L \or\bbb M\or\bbb N\or\bbc O{5} \or
   \bbb P\or\bbc Q{5}\or\bbb R\or\bbc S{4.2}\or\bba T{10.5}\or\bbc U{5}\or
   \bba V{12}\or\bba W{16.5}\or\bba X{11}\or\bba Y{11.7}\or\bba Z{7.5}\fi}}

\def \sysd {\left\{\begin{array}{l}}
\def \sysf {\end{array}\right.}
\def \eqd {\begin{equation}}
\def \eqf {\end{equation}}
\def \f {\end{document}}

\begin{document}

\begin{frontmatter}
\title{Further Remarks on  Strict  Input-to-State \\ Stable Lyapunov
Functions for Time-Varying Systems\thanksref{footnoteinfo}} 

\thanks[footnoteinfo]{
The authors thank Eduardo Sontag for illuminating discussions.
This work was
 supported by Louisiana Board of Regents Contract
LEQSF(2002-04)-ENH-TR-13 (``Visiting Experts Program in
Mathematics''). The work was done in April 2004 while the second
author was a Visiting Expert in the Department of Mathematics at
Louisiana State University (LSU). He thanks LSU for the kind
hospitality he enjoyed during this period. Additional support for
the first author was provided by Louisiana Board of Regents
Contract
 LEQSF(2003-06)-RD-A-12 and NSF Grant 0424011.}

\author{Michael Malisoff}
\address{
Department of Mathematics; 304 Lockett Hall;  Louisiana State
University; Baton Rouge, LA 70803-4918 USA; malisoff@lsu.edu. }
\author{Fr\'ed\'eric Mazenc}
\address{Projet MERE INRIA-INRA; UMR Analyse des Syst\`emes et Biom\'etrie
INRA; 2, pl. Viala; 34060 Montpellier, France;
mazenc@helios.ensam.inra.fr.}


\begin{abstract}
We study the stability properties of a class of time-varying
nonlinear systems. We assume that non-strict input-to-state
stable (ISS) Lyapunov functions for our systems are given and posit
a mild persistency of excitation condition on our given Lyapunov
functions which guarantee the existence of strict ISS Lyapunov functions
for our systems. Next, we provide simple direct constructions
of explicit strict ISS Lyapunov functions for our systems by
applying an integral smoothing method. We illustrate our
constructions using a tracking problem for a rotating rigid body.
\end{abstract}

\end{frontmatter}

\section{Introduction}
The theory of input-to-state stable (ISS) systems plays a central
role in modern non-linear control analysis and controller design
(see \citep{MRS04,MS04,S98,S01,SW95}). The ISS property was
introduced by Sontag in \citep{S89} and an ISS Lyapunov
characterization  was obtained by Sontag and Wang in \citep{SW95}.
The ISS Lyapunov characterization  provides necessary and
sufficient conditions for time-invariant systems to be ISS, in
terms of the existence of so-called strict ISS Lyapunov functions;
see Section \ref{sec2} below for the relevant definitions and
\citep{ELW00} for an extension to time-varying systems. Strict
Lyapunov functions have been used to design stabilizing feedback
laws that render asymptotically controllable systems  ISS to
actuator errors and small observation noise; see \citep{MS04,S01}.
Such control laws are expressed in terms of gradients of Lyapunov
functions and therefore require explicit strict Lyapunov functions
in order to be implemented. This has motivated  a great deal of
research devoted to constructing explicit strict Lyapunov functions.

One obstacle to these constructions is that the known strict
Lyapunov functions from the existence theory are optimal control
value functions, involving a supremum of a cost criterion over
infinitely many possible solution paths (see \citep{BR01,ELW00,SW95,TP00}),
and therefore are not explicit. Although value
functions can often be expressed as unique solutions of
Hamilton-Jacobi (HJ) equations subject to appropriate side
conditions, the usual techniques for computing value functions in
terms of HJ equation solutions can be difficult to implement. For
certain special kinds of systems, strict ISS Lyapunov functions
can be explicitly constructed by ad hoc means. On the other hand,
there are numerous important cases where it is relatively
straightforward to use backstepping or other known methods to
construct explicit {\em non}-strict ISS Lyapunov functions (see our
definitions of ISS and non-strict ISS Lyapunov functions in Section
\ref{sec2}  and Section \ref{sec5} for an explicit example). For
instance, applying the methods of \citep{JN97} to tracking problems
for nonholonomic systems in chained form results in non-strict
Lyapunov functions. The constructions in \citep{MazPrat} also
frequently give rise to non-strict Lyapunov functions.

This motivates the search for techniques for constructing strict
ISS Lyapunov functions for time-varying systems, in terms of known
non-strict ISS Lyapunov functions. This search is the focus of
this note. For time-varying systems with no controls, the paper
\citep{M03} constructed strict globally asymptotically stable
(GAS) Lyapunov functions in terms of given non-strict GAS Lyapunov
functions. Here we  further develop the approach in \citep{M03}.
We provide the necessary background on ISS systems and Lyapunov
functions in Section \ref{sec2}. We then introduce a non-strict
generalization of ISS in which the dissipation rate depends on a
non-negative time-dependent decay parameter. The parameter can be
zero along intervals of positive length. However, when the
parameter is identically one, our non-strict ISS property agrees
with the usual ISS condition. Under a mild non-degeneracy
assumption on this parameter, which is of persistency of
excitation type (see for instance \citep{udel} and \citep{lopa}
for definitions and discussions of the concept of persistency of
excitation), we show that our non-strict ISS property is
equivalent to the existence of a strict ISS Lyapunov function and
is therefore also equivalent to the standard ISS condition. We
prove these equivalences in Section \ref{sec3}. They are proved by
explicitly constructing strict ISS Lyapunov functions. In Section
\ref{sec5}, we illustrate our constructions using a tracking
example. Concluding remarks in Section \ref{secc} end the paper.

\section{Preliminaries}
\label{sec2} Let $\mathcal{K}_\infty$ denote the set of all
continuous functions $\rho:[0,\infty)\to[0,\infty)$ for which (i)
$\rho(0) = 0$ and  (ii) $\rho$ is increasing and unbounded. Let
$\mathcal{K}\mathcal{L}$ denote the set of all continuous
functions $\beta:[0,\infty) \times [0,\infty)\to[0,\infty)$ for
which (1) for each $t\ge 0$, $\beta(\cdot, t)$ is strictly
increasing and $\beta(0,t) = 0$ (2) $\beta(s,\cdot)$ is
non-increasing for each $s\ge 0$, and (3) $\beta(s,t)\to 0$ as
$t\to +\infty$ for each $s\ge 0$.

We  study the stability properties of the fully nonlinear
nonautonomous system
\begin{equation}
\label{sys}
\dot x = f(t,x,u),\; \; \; t\ge 0, \; x\in {\mathbb R}^n,
\; u\in {\mathbb R}^m
\end{equation}
where we  always assume $f$ is  locally Lipschitz  in $(t,x,u)$.
Following \citep{M03}, we also assume $f$ is periodic in $t$,
which means there exists a constant $T>0$ such that
$f(t+T,x,u)=f(t,x,u)$ for all $t\ge 0$, $x\in {\mathbb R}^n$, and
$u\in {\mathbb R}^m$. However, most of our arguments remain valid
if this periodicity assumption is weakened to requiring $f$ to be
uniformly locally bounded in $t$, meaning,
\begin{equation}
\label{boundedness}
\sup\{|f(t,x,u)|: (x,u)\in K, t \ge 0\} < + \infty
\end{equation}
where $|\cdot|$ is the usual Euclidean norm.
The control functions for
our system (\ref{sys}) comprise the set of all measurable locally
essentially bounded functions $\alpha:[0,\infty)\to{\mathbb R}^m$;
we denote this set by $\mathcal{U}$. We let $|\alpha|_I$ denote
the essential supremum of any control $\alpha\in \mathcal{U}$
restricted to any interval $I\subseteq [0,\infty)$. For each
$t_o\ge 0$, $x_o\in {\mathbb R}^n$, and $\alpha\in \mathcal{U}$,
we let $I\ni t\mapsto \phi(t;x_o, t_o, \alpha)$ denote the unique
trajectory of (\ref{sys}) for the input $\alpha$ satisfying
$x(t_o) = x_o$  and defined on its maximal interval $I\subseteq [t_o,\infty)$.
This trajectory will be denoted by $\phi$ when this would not lead to confusion.
We say that $f$ is {\em forward complete} provided each such trajectory $\phi$ is
defined on all of $[t_o,\infty)$.

A $C^1$ function $V:[0,\infty) \times {\mathbb R}^n \to [0,\infty)$
is said to be of class ${\rm UPPD}$ (written $V \in {\rm UPPD}$)
provided it is uniformly proper and positive definite, which means
there exist $\alpha_1, \alpha_2, \alpha_3 \in \mathcal{K}_\infty$ such
that, for all $t \ge 0, x \in {\mathbb R}^n$,
\begin{equation}
\label{uppd} \alpha_1(|x|) \le V(t,x)\le \alpha_2(|x|), \; \;
|\nabla V(t,x)|\le \alpha_3(|x|).
\end{equation}
We say that $V$ has period $\tau$ in $t$ provided there exists a
constant $\tau>0$ such that $V(t+\tau,x)=V(t,x)$ for all $t\ge 0$
and $x\in {\mathbb R}^n$; in this case, the bound on $\nabla V$ in
(\ref{uppd})  is redundant.  We  assume $\alpha_1$ and $\alpha_2$
in (\ref{uppd}) are $C^1$, e.g., by taking
$\alpha_2(s)=\int_o^{\scriptscriptstyle s}\alpha_3(r)dr$ and
minorizing $\alpha_1$ by a $C^{\scriptscriptstyle 1}$ function of
class $\mathcal{K}_\infty$.
 Given $V\in {\rm UPPD}$, we set
\[\dot V(t,x,u):=\frac{\partial V}{\partial t}(t,x)+
\frac{\partial V}{\partial x}(t,x) f(t,x,u).\] Notice that
$s\mapsto\sup\{|\dot V(t,x,u)|:  t\ge 0, |x|\le \chi(s),|u|\le s\}
+s$ is of class $\mathcal{K}_\infty$ for each $\chi\in
\mathcal{K}_\infty$ (by (\ref{boundedness})-(\ref{uppd})). We let
$\mathcal{P}$ denote the set of all continuous functions
$p:{\mathbb R} \to [0,\infty)$ that admit constants
$\tau,\varepsilon,\bar p>0$ for which
\begin{equation}
\label{calp}
\begin{array}{rcl}
\int_{t-\tau}^tp(s)ds\; \ge \; \varepsilon\; {\rm \and\ }  \; p(t)\le \bar p\; , \; \; \forall t\ge 0.
\end{array}
\end{equation}
We write $p\in \mathcal{P}(\tau,\varepsilon,\bar p)$ to indicate that
(i) $p\in \mathcal{P}$ and (ii) $\tau,\varepsilon,\bar p>0$ are
constants such that (\ref{calp}) holds.  In particular, any
continuous periodic function $p:{\mathbb R}\to[0,\infty)$ that  is
not identically zero admits constants $\tau,\varepsilon,\bar p>0$
satisfying (\ref{calp}).  On the other hand, (\ref{calp}) also
allows non-periodic $p$ with arbitrarily large null sets, e.g., for
fixed $r > 0$, set $p_r(t) = (1 + e^{-t})\max\{0, \sin^3(\frac{t}{r})\}$.
The elements of $\mathcal{P}$
serve as the decay rates for our non-strict Lyapunov functions as
follows:
\begin{defn}
Let $p\in \mathcal{P}$. A function $V\in {\rm UPPD}$ is called
an {\rm ISS(p) Lyapunov function} for (\ref{sys}), provided there exist
$\chi\in \mathcal{K}_\infty$ and  $\mu\in \mathcal{K}_\infty\cap C^1$
such that
\begin{equation}
\label{nonstr} |x| \ge \chi(|u|) \Rightarrow \dot V(t,x,u) \le -
p(t)\mu(|x|)\;  \; \; \forall t\ge 0.
\end{equation}
An ISS(p) Lyapunov function for (\ref{sys}) and $p(t)\equiv 1$
is also called a {\rm strict  ISS Lyapunov function}.
\end{defn}
Notice that (\ref{nonstr}) allows $\dot V(t,x,u)=0$ for those $t$
where $p(t) = 0$. This corresponds to allowing $V$  to non-strictly
decrease along the solutions $\phi$ of (\ref{sys}).
\begin{defn}
Let $p\in \mathcal{P}$.  We say that (\ref{sys}) is  {\rm ISS(p)},
or that it is {\rm input-to-state stable (ISS) with decay rate}
$p$, provided there exist $\beta\in \mathcal{K}\mathcal{L}$ and
$\gamma \in \mathcal{K}_\infty$ such that for all $t_o\ge 0$,
$x_o\in {\mathbb R}^n$, $u_o\in \mathcal{U}$ and $h\ge 0$,
\begin{equation}
\label{siss}
\begin{array}{rcl}
|\phi(t_o+h; x_o, t_o, u_o)| & \le & \beta\left(|x_o|, \int_{t_o}^{t_o+h} p(s)ds\right)
\\
& & + \gamma\left(|u_o|_{[t_o,t_o+h]}\right).
\end{array}
\end{equation}
If (\ref{sys}) is ISS(p) with  $p \equiv 1$, then we say that (\ref{sys}) is {\rm ISS}.
\end{defn}
Notice that ISS(p) systems are automatically forward complete. We
also study dissipation-type decay conditions as follows:
\begin{defn}
\label{zpb} Let $p\in \mathcal{P}$. A function $V\in {\rm UPPD}$
is called a {\rm non-strict dissipative Lyapunov function} for
(\ref{sys}) and $p$, or a {\rm DIS(p) Lyapunov function}, provided
there exist $\Omega\in \mathcal{K}_\infty$ and $\mu\in
\mathcal{K}_\infty\cap C^1$ such that, for all $t \ge 0, x \in
{\mathbb R}^n, u \in {\mathbb R}^m$
\begin{equation}
\label{i3}\dot V(t,x,u) \; \le - p(t)\mu(|x|) + \Omega(|u|) \; \; .
\end{equation}
A DIS(p) Lyapunov function for (\ref{sys}) and $p(t)\equiv 1$ is
also called a {\rm strict  DIS Lyapunov function}.
\end{defn}
\begin{rem}
\label{rte} Definition \ref{zpb} is a nonlinear version of the
property used in \citep{lopa} to ensure the global uniform
exponential stability of time-varying linear systems belonging to
a specific family of systems. Thus, the explicit construction of a
strict DIS Lyapunov function in terms of a given DIS(p) Lyapunov
function we present in the next section, extends \citep{lopa}
where only linear systems are studied and no strict Lyapunov
function is constructed.
\end{rem}

We use the following elementary observations:
\begin{lem}
\label{plem}
Let $\tau,\varepsilon, \bar p > 0$ be constants and
$p\in \mathcal{P}(\tau,\varepsilon,\bar p)$ be given.  Then: \newline
\noindent
(i) $0 \leq \int_{t-\tau}^t\left(\int_s^tp(r)dr\right)ds \le \frac{\tau^2\bar p}{2}$
for all $t\ge 0$ and\newline
\noindent
(ii) $[0,\infty)\ni h\mapsto {\underline
p}(h)=\inf\left\{\int_t^{t+h}p(r)dr: t\ge 0\right\}$ is
continuous, non-decreasing, and unbounded.
\end{lem}
We leave the proof of this lemma to the reader as a simple
exercise.
\section{Equivalent Characterizations of Non-Strict ISS}
\label{sec3}
We next relate the Lyapunov functions and stability notions we
introduced in the last section. We show that ISS(p)
 is equivalent to the existence of an ISS(p) Lyapunov function
 and the existence of a strict ISS Lyapunov function. Our proof
 explicitly constructs a  strict ISS
Lyapunov function for (\ref{sys})  in terms of a given DIS(p)
Lyapunov function.  Moreover, if $p\in \mathcal{P}(\tau,\varepsilon,
\bar p)$ and our given DIS(p) Lyapunov function both have period
$\tau$, then the strict ISS Lyapunov function we construct also
has period $\tau$. We next prove:
\begin{thm}
\label{mainthm}
Let $p\in \mathcal{P}$ and $f$ be as above. The following are equivalent: \newline
\noindent
$(C_1)$ $f$ admits an ISS(p) Lyapunov function.\newline
\noindent
$(C_2)$ $f$ admits a strict ISS Lyapunov function.\newline
\noindent
$(C_3)$ $f$ admits a DIS(p) Lyapunov function.\newline
\noindent
$(C_4)$ $f$ admits a strict DIS Lyapunov function.\newline
\noindent
$(C_5)$ $f$ is ISS(p). \newline
\noindent
$(C_6)$ $f$ is ISS.
\end{thm}
\noindent  We prove the following implications: $(C_1)
\Rightarrow (C_2) \Rightarrow (C_4)\Rightarrow (C_1)$,
$(C_3)\Leftrightarrow(C_4)$, $(C_2)\Leftrightarrow(C_6)$, and
$(C_5)\Leftrightarrow(C_6)$.
We fix  $\tau,\varepsilon, \bar p>0$  such that $p\in
\mathcal{P}(\tau,\varepsilon,\bar p)$.\newline \noindent {\bf Step
1:} $(C_1) \Rightarrow (C_2)$. If $(C_1)$ holds, then we can find
an ISS(p) Lyapunov function $V$ for $f$, and therefore
$\alpha_{1}, \alpha_{2}\in \mathcal{K}_{\infty}\cap C^1$
satisfying (\ref{uppd}) and $\chi\in \mathcal{K}_\infty$ and $\mu
\in \mathcal{K}_{\infty}\cap C^1$ satisfying (\ref{nonstr}). Set
\begin{equation}
\label{choices}
\begin{array}{rcl}
\tilde\alpha_2(s) & := & \max\left\{\frac{\tau\bar
p}{2},1\right\}(\alpha_2(s)+\mu(s)+s),
\\
w(s) & := & \frac{1}{4\tau}\mu(\tilde \alpha^{-1}_2(s)).
\end{array}
\end{equation}
Then $\tilde \alpha_2,\tilde \alpha^{-1}_2 \in
\mathcal{K}_\infty\cap C^1$. Since $V(t,x)\le \tilde
\alpha_2(|x|)$ for all $t\ge 0$ and $x\in {\mathbb R}^n$, the
following holds for all $t\ge 0$:
\begin{equation}
\label{g1}  |x| \ge \chi(|u|)  \Rightarrow  \dot V(t,x,u) \le -
p(t)\mu(\tilde \alpha_2^{-1}(V(t,x))) .
\end{equation}
Note too that $w\in \mathcal{K}_{\infty}\cap C^1$. We later use
the fact that
\begin{equation}
\label{ig7} \begin{array}{lll} \displaystyle 0\; \le\;
w'(s)&\le&\displaystyle \frac{\mu'(\tilde\alpha^{-1}_2(s))}{4\tau
\max\{\frac{\tau\bar p}{2},1\}(\mu'(\tilde \alpha^{-1}_2(s))+1)}\\
&\le& \displaystyle\frac{1}{2\tau^2\bar p}\end{array}
\end{equation}
for all $s\ge 0$. Consider the UPPD function
\begin{equation}
\label{vs}
V^\sharp(t,x) = V(t,x) + \xi(t) w(V(t,x))
\end{equation}
with $\xi(t) = \int_{t-\tau}^{t}\left(\int_s^t p(r)\, dr\,\right)
ds$. Then
\[
\begin{array}{rcl}
\dot{V}^\sharp(t,x,u) & = & [1 + \xi(t) w'(V(t,x))]\dot{V}(t,x,u)
\\
& &  + \left[\tau p(t) - \int_{t - \tau}^{t} p(r)\, dr\right]
w(V(t,x))
\end{array}\]
follows from a simple calculation. When $|x| \ge \chi(|u|)$,
condition (\ref{g1}) gives ${\scriptstyle \dot V}(t,x,u)\le 0$ and
therefore also
\[
\begin{array}{rcl}
\dot{V}^\sharp(t,x,u) & \leq & - p(t)\mu(\tilde
\alpha^{-1}_2(V(t,x)))
\\
& +&  \left[\tau p(t) - \int_{t - \tau}^{t} p(r)\, dr\right]
\frac{1}{4\tau}\mu(\tilde \alpha_2^{-1}(V(t,x)))
\\
& \leq & - \frac{3}{4}p(t)\mu(\tilde \alpha^{-1}_2(V(t,x)))
\\
& -&  \left(\int_{t - \tau}^{t} p(r)\, dr\right)
\frac{1}{4\tau}\mu(\tilde \alpha_2^{-1}(V(t,x)))
\\
& \leq & -
\frac{\varepsilon}{4\tau}\mu(\tilde\alpha_2^{-1}(\alpha_1(|x|)))\;
\; \forall t\ge 0.
\end{array}
\]
Since $\mu\circ\tilde \alpha_2^{-1}\circ\alpha_1\in C^1\cap
\mathcal{K}_\infty$, it follows that $V^\sharp$ is a strict ISS
Lyapunov function for (\ref{sys}).\newline \noindent {\bf Step 2:}
$(C_2) \Rightarrow (C_4)$. Assume $(C_2)$, so $f$ admits a strict
ISS Lyapunov function $V$. Let $\mu$ and $\chi$ satisfy  condition
(\ref{nonstr}) with $p\equiv 1$. Then the strict dissipative
condition (\ref{i3}) with $p\equiv 1$ follows by choosing any
$\Omega \in \mathcal{K}_\infty$ satisfying
\[
\Omega(s) \ge \displaystyle\max_{\{t\ge 0, |x|\le \chi(s), |u|\le
s\}}\{\dot V(t,x,u) + \mu(|x|)\} \; \; \forall s\ge 0.
\]
Such an $\Omega$ exists by our assumptions
(\ref{boundedness})-(\ref{uppd}). Therefore, $V$ is itself a
strict DIS Lyapunov function for $f$. \newline \noindent {\bf Step
3:} $(C_4) \Rightarrow (C_1)$. Assume $(C_4)$, so $f$ admits a
strict  DIS Lyapunov function $V$. Let $\mu,\Omega\in
\mathcal{K}_\infty$ satisfy (\ref{i3}) with $p\equiv 1$; then if
$|x| \geq \chi(|u|):=\mu^{-1}(2\Omega(|u|))$, then
\[\dot V(t,x,u) \le - \frac{1}{2} \mu(|x|),\; \; {\rm so}\; \;
\dot V(t,x,u) \le - \frac{p(t)}{2\bar p} \mu(|x|)\] for all $t\ge
0$.
Therefore, $V$ is also an ISS(p) Lyapunov function for $f$, so
$(C_1)$ is satisfied.\newline \noindent {\bf Step 4:}
$(C_3)\Leftrightarrow(C_4)$. Since $p\in \mathcal{P}$ is bounded,
we easily conclude that $(C_4)$ implies $(C_3)$. Conversely,
assume $V\in {\rm UPPD}$ is a DIS(p) Lyapunov function for $f$ and
$\alpha_1,\alpha_2,\mu,\Omega\in \mathcal{K}_\infty$ satisfy
(\ref{uppd}) and the DIS(p) requirements. Define $\tilde \alpha_2,
w\in \mathcal{K}_\infty\cap C^1$ and $V^\sharp$ by (\ref{choices})
and (\ref{vs}).   As before, when $\tilde \mu=\mu\circ \tilde
\alpha^{\scriptscriptstyle -1}_2$, we have ${\scriptstyle \dot
V}(t,x,u) \le - p(t)\tilde \mu(V(t,x)) + \Omega(|u|)$  for all $t
\ge 0, x \in {\mathbb R}^n, u \in {\mathbb R}^m$. It follows from
Lemma \ref{plem}(i) and (\ref{ig7}) that
\begin{equation}\label{bounds}
1 + \xi(t) w'(V(t,x)) \in \left[1,\frac{5}{4}\right]\; , \; \;
\forall t\ge 0, x\in {\mathbb R}^n.
\end{equation}
Since  $w = \frac{1}{4\tau}\tilde \mu$, we deduce that
\begin{equation}
\begin{array}{rcl}
\dot{V}^\sharp & \leq & - p(t)\tilde{\mu}(V(t,x)) +
\frac{5}{4}\Omega(|u|)\nonumber
\\
& +&  \tau p(t) w(V(t,x))\nonumber - \left(\int_{t - \tau}^{t}
p(r)dr\right) w(V(t,x))\nonumber
\\
& \leq & - \varepsilon w(\alpha_1(|x|)) + \frac{5}{4}\Omega(|u|). \nonumber
\end{array}
\end{equation}
Since $w\circ\alpha_1\in C^1\cap \mathcal{K}_{\infty}$, it follows
that $V^\sharp$ is the desired strict DIS Lyapunov
function.\newline \noindent {\bf Step 5:} $(C_2)\Leftrightarrow
(C_6)$. The implication $(C_2)\Rightarrow (C_6)$ follows from
\citep[Theorem 4.19, p.176]{K02}.  (In \citep{K02}, the controls
are bounded piecewise continuous functions $\alpha:[0,\infty)\to
{\mathbb R}^m$, but the result from \citep{K02} can be extended to
our general control set $\mathcal{U}$ using a standard denseness
argument (see e.g. Remark C.1.2 and the proof of Theorem 1 in
\citep{S98a}).) The converse was announced in \citep[Theorem
1]{ELW00} and can be deduced from \citep{BR01} as follows. If $f$
is ISS, then \citep{SW95} provides $\chi\in \mathcal{K}_\infty$
such that the constrained input system $\dot x =
f_\chi(t,x,d):=f(t,x,d\chi^{-1}(|x|))$, $|d|\le 1$ is uniformly
globally asymptotically stable (UGAS); i.e., there exists $\beta
\in \mathcal{K}\mathcal{L}$ such that for each $t_o\ge 0$ and
$x_o\in {\mathbb R}^n$ and  each trajectory $y$ of $f_\chi$
satisfying $y(t_o)=x_o$, we have $|y(t_o+h)|\le \beta(|x_o|,h)$
for all $h\ge 0$.  By minorizing $\chi^{-1}$, we can assume it is
$C^1$. This means the locally Lipschitz set-valued dynamics
$F(t,x)=\{f(t,x,u): \chi(|u|)\le |x|\}$ is UGAS, as is its
convexification $\overline{\rm co}(F)$, namely $(t,x)\mapsto
\overline{\rm co}\{F(t,x)\}$ where $\overline{\rm co}$ denotes the
closed convex hull (cf. \citep[Proposition 4.2]{BR01}). Since
$\overline{\rm co}(F)$ is continuous and compact and convex
valued, and since we are assuming $f$ is periodic in $t$,
\citep[Theorem 4.5]{BR01} provides a time-periodic $V\in {\rm
UPPD}$ such that, for all $x\in {\mathbb R}^n, \; t\ge 0,\; w \in
F(t,x)$,
\[
\frac{d}{dt}V(t,x)+\frac{d}{dx} V(t,x)w\le - V(t,x).
\]
Recalling the definition of $F$ and assuming (without loss of
generality) that $V$ satisfies (\ref{uppd}) with $\alpha_1\in
\mathcal{K}_\infty\cap C^1$,
\[
\begin{array}{l}
|x|\ge \chi(|u|)\; \Rightarrow\; f(t,x,u)\in F(t,x)\\
\Rightarrow\; \dot V(t,x,u)\le -V(t,x)\le
-\alpha_1(|x|)\end{array}\] for all $t\ge 0$, so $V$ is the
desired strict ISS Lyapunov function for $f$. This establishes
$(C_6)\Rightarrow (C_2)$.\newline \noindent {\bf Step 6:}
$(C_5)\Leftrightarrow (C_6)$. Assuming $(C_6)$, there are
$\beta\in \mathcal{K}\mathcal{L}$  such that for all $t_o\ge 0$,
$x_o\in {\mathbb R}^n$, $u_o\in \mathcal{U}$, and $h\ge 0$,
\[
\begin{array}{lll}
|\phi(t_o+h; x_o, t_o, u_o)|& \le & \beta(|x_o|, \bar p h)
+\gamma(|u_o|_{[t_o,t_o+h]})
\\
& \le & \beta(|x_o|,\int_{t_o}^{t_o+h} p(s)ds ) \\&+&
\gamma(|u_o|_{[t_o,t_o+h]}),\end{array}
\]
where $\phi$ is the trajectory of $f$ we defined
in Section \ref{sec2}. Therefore, $f$ is ISS(p) so $(C_6)\Rightarrow (C_5)$.
Conversely, if $f$ is ISS(p), then we can find
$\beta\in \mathcal{K}\mathcal{L}$ such that for all
$t_o \ge 0$, $x_o \in {\mathbb R}^n$, $u_o \in \mathcal{U}$, and  $h \ge 0$,
\[
\begin{array}{lll}
|\phi(t_o+h; x_o, t_o, u_o)|& \le & \beta\left( |x_o|,
\int_{t_o}^{t_o+h} p(s)ds\right)
\\
& & + \gamma(|u_o|_{[t_o,t_o+h]})
\\
& \le & \beta\left(|x_o|,\underline{p}(h)\right)+\gamma(|u_o|_{[t_o,t_o+h]}).
\end{array}
\]
By Lemma \ref{plem}(ii) ,  $\hat \beta(s,t):=
\beta(s,\underline{p}(t))\in \mathcal{K}\mathcal{L}$, so $(C_5)\Rightarrow
(C_6)$, as desired. This proves Theorem \ref{mainthm}.

\medskip

\begin{rem}
Observe that if the functions $V$, $\alpha_2$, $\mu$, $p$ are of
class $C^k$, where $k$ is a positive integer or $\infty$, then the
particular function $\tilde{\alpha}_2$ in (\ref{choices}) we have
chosen implies that the function $V^\sharp(t,x)$ is of class
$C^k$.
\end{rem}
\begin{rem}
Our proof of Theorem \ref{mainthm} shows that if $V$ is a strict
ISS Lyapunov function for $f$, then $V$ is also a strict DIS
Lyapunov function for $f$. The preceding implication is no longer
true if our boundedness requirement (\ref{boundedness})  on $f$ is
dropped, as illustrated by the following example from
\citep{ELW00}: Take the one-dimensional single input system $\dot
x=f(t,x,u):=-x+(1+t)q(u-|x|)$, where $q:{\mathbb R}\to{\mathbb R}$
is any $C^1$ function for which $q(r)\equiv  0$ for $r\le 0$ and
$q(r)>0$ otherwise.
Then $V(x)=x^2$ is a strict ISS Lyapunov function for the system
since $|x|\ge |u| \Rightarrow {\scriptstyle \dot V}\le -x^2$ but
$V$ does
 not satisfy the strict DIS condition (\ref{i3}) for any
choices of $\mu$ and $\Omega$.  This does not contradict our
results because (\ref{boundedness}) is not satisfied. This
contrasts with the time-invariant case where strict ISS Lyapunov
functions are automatically strict DIS Lyapunov functions.
\end{rem}

\section{Illustration}
\label{sec5} We next use our results to construct a strict ISS
Lyapunov function for a tracking problem for a rotating rigid body
(see \citep{C84,SSPJ95,SS97} for the background and motivation for
this problem). Following Lefeber \citep[p.31]{L00}, we only
consider the dynamics of the velocities, which, after a change of
feedback, are
\begin{equation}
\label{transformed}
\dot{\omega}_1 = \delta_1 + u_1\; ,\; \;
\dot{\omega}_2 = \delta_2 + u_2\; ,\; \; \dot{\omega}_3 = \omega_1\omega_2.
\end{equation}
where $\delta_1$ and $\delta_2$ are the inputs and $u_1$ and $u_2$ are
the disturbances.
We consider the reference state trajectory
\begin{equation}
\label{zl1}
\omega_{1r}(t) = \sin(t) \; ,\; \; \omega_{2r}(t)= \omega_{3r}(t) = 0
\end{equation}
but our method applies to more general reference trajectories as
well; see Remark \ref{lasr} below. The substitution
$\tilde{\omega}_i(t) = \omega_i(t) - \omega_{ir}(t)$ transforms
(\ref{transformed}) into the error equations
\begin{equation}
\label{cnh}
\begin{array}{rcl}
\dot{\tilde{\omega}}_1 & = & \delta_1 + u_1 - \cos(t) \; ,
\\
\dot{\tilde{\omega}}_2 & = & \delta_2 + u_2 \; ,
\\
\dot{\tilde{\omega}}_3 & = & (\tilde{\omega}_1 + \sin(t))\tilde{\omega}_2\; .
\end{array}
\end{equation}
By applying the backstepping approach as it is applied in
\citep{JN97}, or through direct calculations, one  shows that the
derivative of the class UPPD function
\begin{equation}
\label{tjb}
V(t,\tilde{\omega}) = \frac{1}{2}\left[\tilde{\omega}_1^2
+ (\tilde{\omega}_2 + \sin(t)\tilde{\omega}_3)^2 + \tilde{\omega}_3^2\right]
\end{equation}
with $\tilde{\omega} = (\tilde{\omega}_1, \tilde{\omega}_2, \tilde{\omega}_3)^\top$
along the trajectories of (\ref{cnh}) in closed-loop with the control laws
\begin{equation}
\label{poc}
\begin{array}{rcl}
\delta_1(t,\tilde{\omega}) & = & - \tilde{\omega}_1 - \tilde{\omega}_2\tilde{\omega}_3 + \cos(t)
\\
\delta_2(t,\tilde{\omega}) & = & - [1 + \sin(t)\tilde{\omega}_1 + \sin^2(t)]\tilde{\omega}_2
\\
& & - (2\sin(t) + \cos(t))\tilde{\omega}_3
\end{array}
\end{equation}
satisfies
\begin{equation}
\label{ufg}
\begin{array}{rcl}
\dot{V} & = & - \tilde{\omega}_1^2 - (\tilde{\omega}_2 + \sin(t)\tilde{\omega}_3)^2
- \sin^2(t)\tilde{\omega}_3^2
\\
& & + \tilde{\omega}_1 u_1 + (\tilde{\omega}_2 + \sin(t)\tilde{\omega}_3) u_2
\\
& \leq & - \frac{1}{2}\tilde{\omega}_1^2 - \frac{1}{2}(\tilde{\omega}_2 + \sin(t)\tilde{\omega}_3)^2
- \sin^2(t)\tilde{\omega}_3^2
\\
& & + \frac{1}{2}(u_1^2 + u_2^2)
\\
& \leq & - p(t)\tilde\mu(V(\tilde{\omega}))
+ \Omega(|u|)
\end{array}
\end{equation}
with $u = (u_1,u_2)^\top \in {\mathbb R}^2$, $p(t) = \sin^2(t)$,
$\tilde\mu(s) = s$ and $\Omega(s) = \frac{1}{2}s^2$. Therefore $V$
is a DIS(p) Lyapunov function for (\ref{cnh}) in closed-loop with
the control laws (\ref{poc}). Observe that, in this case, $p\in
\mathcal{P}(\pi,\pi/2,1)$. Setting $\tau = \pi$ and $w(s) =
\frac{1}{8\tau}\tilde\mu(s) = \frac{s}{8\pi}$, it follows that
(\ref{bounds}) also holds. Therefore, Steps 3-4 from our proof of
Theorem \ref{mainthm} show
\[\begin{array}{ccc}
\! \! & V^\sharp(t,\tilde{\omega}) = V(t,\tilde{\omega}) +
\left[\int_{t-\tau}^t\left(\int_s^t p(r) dr\right) ds\right]
w(V(t,\tilde{\omega}))
\\
& = \left[1 + \frac{\pi}{32} - \frac{1}{32}
\sin(2t)\right]V(t,\tilde{\omega})  
\end{array}
\]
is a strict DIS Lyapunov function and also a strict ISS Lyapunov
function for the system (\ref{cnh}) in closed-loop with
the control laws (\ref{poc}). \medbreak\medbreak \noindent

\begin{rem}
\label{lasr} We chose to work with the reference trajectory
(\ref{zl1}) because it leads to the simple error equations
(\ref{cnh}). However, one can easily check that a strict ISS
Lyapunov function can be constructed for any reference state
trajectory $(\omega_{1r}(t),\omega_{2r}(t), \omega_{3r}(t))$ such
that \[\begin{array}{l}\sup_t
\left|\int_{0}^{t}\omega_{1r}(s)\omega_{2r}(s) ds\right| <
\infty\; \; {\rm and}\\  \int_{t - \tau}^{t}[\omega^2_{1r}(s) +
\omega^2_{2r}(s)] ds \geq \varepsilon \; , \; \forall t\ge
\tau\end{array}\]
for some constants $\tau,\varepsilon>0$.
\end{rem}

\section{Conclusion}
\label{secc}

For ISS time-varying systems, we provided  explicit strict
Lyapunov function {\em constructions} that can easily be performed
in practice. The knowledge of these Lyapunov functions allows us
to extend the well-known and useful theory of ISS systems to a
broad class of  time-varying nonlinear dynamics.  We conjecture
that a discrete-time version of our main result can be proved.

\thebibliography{xx}

\harvarditem[Bacciotti \& Rosier]{Bacciotti \& Rosier}{2001}{BR01}
Bacciotti, A. \& Rosier, L. (2001). {\em Liapunov Functions and
Stability in Control Theory.} Lecture Notes in Control and Inform.
Sci.  Vol. 267, Springer-Verlag London, Ltd., London.

\harvarditem[Crouch]{Crouch}{1984}{C84} Crouch, P. (1984).
Spacecraft attitude control and stabilization: applications of
geometric control theory to rigid body models. {\em IEEE Trans.
Automat. Control}, 29(4), 321-331.

\harvarditem[Edwards {\em et al.}]{Edwards, Lin \&
Wang}{2000}{ELW00} Edwards, H., Lin, Y. \& Wang, Y. (2000). On
input-to-state stability for time-varying nonlinear systems. {\em
Proceedings of the 39th IEEE Conference on Decision and Control,
Sydney, Australia.}

\harvarditem[Jiang \& Nijmeijer]{Jiang \& Nijmeijer}{1997}{JN97}
Jiang, Z.-P. \& Nijmeijer, H. (1997). Tracking control of mobile
robots: a case study in backstepping. {\em Automatica}, 33(7),
1393-1399.

\harvarditem[Khalil]{Khalil}{2002}{K02}
Khalil, H. (2002).
{\em Nonlinear Systems, Third Edition.} Prentice Hall, 2002.

\harvarditem[Lefeber]{Lefeber}{2000}{L00}
Lefeber, E., (2000).
{\em Tracking Control of Nonlinear Mechanical Systems.}
PhD Thesis, University of Twente, Enschede, The Netherlands, April 2000.
(On-line at $\mathtt{http://se.wtb.tue.nl/\sim lefeber/}$.)

\harvarditem[Loria \& Panteley]{Loria \& Panteley}{2002}{lopa}
Loria, A. \& Panteley, E. (2002). Uniform exponential stability of
linear time-varying systems: revisited. {\em Systems \& Control
Letters}, 47(1), 13-24.

\harvarditem[Loria {\em et al.}]{Loria, Panteley, Popovic \&
Teel}{2002}{udel} Loria, A., Panteley, E., Popovic, D. \& Teel, A.
(2002). $\delta$-persistency of excitation: a necessary and
sufficient condition for uniform attractivity. {\em Proceedings of
the 41st IEEE Conference on Decision and Control, Las Vegas, NV.}

\harvarditem[Malisoff {\em et al.}]{Malisoff, Rifford \&
Sontag}{2004}{MRS04} Malisoff, M., Rifford L. \& Sontag E. (2004).
Global asymptotic controllability implies input-to-state
stabilization. {\em SIAM Journal on Control and Optimization,}
42(6), 2221-2238.

\harvarditem[Malisoff \& Sontag]{Malisoff, Rifford \&
Sontag}{2004}{MS04} Malisoff, M. \& Sontag, E. (2004). {\em
Asymptotic controllability and input-to-state stabilization:  The
effect of actuator errors.} In: Optimal Control, Stabilization,
and Nonsmooth Analysis,  Lecture Notes in Control and Inform.
Sci., Springer-Verlag, Heidelberg.  Vol. 301, 155-171.

\harvarditem[Mazenc]{Mazenc}{2003}{M03} Mazenc, F. (2003). Strict
Lyapunov functions for time-varying systems. {\em Automatica},
39(2), 349-353.

\harvarditem[Mazenc \& Praly]{Mazenc \& Praly}{2000}{MazPrat}
Mazenc, F. \& Praly, L. (2000). Asymptotic Tracking of a State
Reference for Systems with a Feedforward Structure. {\em
Automatica}, 36(2), 179-187.

\harvarditem[Morin \& Samson]{Morin \& Samson}{1997}{SS97} Morin,
P. \& Samson, C. (1997). Time-varying exponential stabilization of
a rigid spacecraft with two controls. {\em IEEE Trans. Automatic
Control}, 42(4), 528-534.

\harvarditem[Morin {\em et al.}]{Morin, Samson, Pomet \&
Jiang}{1995}{SSPJ95} Morin, P., Samson, C., Pomet, J.-B. \& Jiang,
Z.-P. (1995). Time-varying feedback stabilization of the attitude
of a rigid spacecraft with two controls. {\em Systems \& Control
Letters}, 25(5), 375-385.

\harvarditem[Sontag]{Sontag}{1989}{S89} Sontag, E. (1989). Smooth
stabilization implies coprime factorization. {\em IEEE Trans.
Automatic Control}, 34(4), 435-443.

\harvarditem[Sontag]{Sontag}{1998}{S98} Sontag, E. (1998).
Comments on integral variants of ISS. {\em Systems \& Control
Letters}, 34(1-2), 93-100.

\harvarditem[Sontag]{Sontag}{1998}{S98a}
Sontag, E. (1998).
{\em Mathematical Control Theory. Deterministic
Finite-Dimensional Systems. Second Edition.} Texts in Applied
Mathematics 6. Springer-Verlag, New York, 1998.

\harvarditem[Sontag]{Sontag}{2001}{S01} Sontag, E. (2001). {\em
The ISS philosophy as a unifying framework for stability-like
behavior}.  In: Nonlinear Control in the Year 2000, Vol. 2,
Lecture Notes in Control and Inform. Sci.,  Springer, London. Vol.
259, 443--467.

\harvarditem[Sontag \& Wang]{Sontag \& Wang}{1995}{SW95} Sontag,
E. \& Wang, Y. (1995). On characterizations of the input-to-state
stability property. {\em Systems \& Control Letters}, 24(5),
351-359.

\harvarditem[Teel \& Praly]{Teel \& Praly}{2000}{TP00} Teel, A. \&
Praly, L. (2000). A smooth Lyapunov function from a
class-$\mathcal{K}\mathcal{L}$ estimate involving two positive
semidefinite functions. {\em ESAIM Control Optim. Calc. Var.}, 5,
313-367.

\endthebibliography

\end{document}